\documentclass[a4paper,11pt]{amsart}

\usepackage{amsmath, amssymb, amsthm}
\usepackage[top=3cm, bottom=3cm, left=2.5cm, right=2.5cm]{geometry}
\usepackage[colorlinks=true, linkcolor=blue, citecolor = blue]{hyperref}

\newtheorem{thm}{Theorem}[section]
\newtheorem{assumpt}{Assumption}[section]
\newtheorem{rmk}{Remark}[section]

\newcommand{\eps}{\varepsilon}
\newcommand{\expt}{\mathbb{E}}

\newcommand{\R}{\mathbb{R}}
\newcommand{\C}{\mathcal{C}}
\newcommand{\Ind}{\mathbf{1}_{\{0\}}}
\newcommand{\la}{\langle}
\newcommand{\ra}{\rangle}
\newcommand{\Psiz}{\hat{\Psi}}
\newcommand{\pj}{\Pi}
\newcommand{\zed}{\mathfrak{z}^\eps}
\newcommand{\zedx}{\mathfrak{z}}
\newcommand{\gnoise}{\xi}

\newcommand{\ubar}[1]{\underline{#1}}

\title[LDP for DDE at an instability]{Large deviations in presence of small noise for delay differential equations at an instability }
\author{Nishanth Lingala}

\begin{document}

\begin{abstract}
We consider delay differential equations (DDE) that are on the verge of an instability, i.e. the characteristic equation for the linearized equation has one root as zero and all other roots have negative real parts. In presence of small mean-zero noise, we study the large deviations from the corresponding deterministic system. Using spectral theory for DDE it is easy to see that, the projection on to the one dimensional space corresponding to the zero root is exponentially equivalent with the original process. For the one-dimensional process we make the observation that the results of Freidlin-Wentzell apply.
\end{abstract}

\maketitle


\section{Introduction}
We consider $\R^n$ valued processes governed by delay differential equations (DDE) of the form
\begin{align}\label{eq:TB:ldp:main}
\dot{x}(t)=L_0(\pj_tx)+\eps G(\pj_t x)+\eps F(\pj_t x)\sigma(\gnoise_t),
\end{align}
where 
\begin{itemize}
\item $\pj_t$ is the \emph{segment extractor} defined by $(\pj_tf)(\theta)=f(t+\theta)$ for $\theta\in [-r,0]$ where $r>0$ is the maximal delay in the system; note that
$\pj_t:C([-r,\infty);\R^n)\to \C:=C([-r,0];\R^n)$
\item $L_0,G,F:\C\to\R^n$, with $L_0$ being linear, and $G,F$ being bounded with bounded derivatives
\item $\sigma$ is a bounded mean zero $\R$-valued function of the Markov noise $\gnoise$
\item $\eps\ll 1$ is a small number signifying a perturbation.
\end{itemize}
We assume that there exists a bounded matrix-valued function $\mu:[-r,0]\to \R^{n\times n}$, continuous from the left on the interval $(-r,0)$ and normalized with $\mu(0)=0_{n\times n}$, such that 
\begin{align}\label{eq:TB:ldp:L0matrixrep}
L_0\eta = \int_{[-r,0]}d\mu(\theta)\eta(\theta), \quad \forall \eta \in \C.
\end{align}
This is not a restriction: every continuous linear operator $L_0$ has such a representation. 

We make the following assumption on $L_0$ to reflect an instability scenario:  
\begin{assumpt}\label{ass:TB:ldp:assumptondetsys}
Define $$\Delta(\lambda)\,=\,\lambda I_{n\times n}-\int_{[-r,0]}d\mu(\theta)e^{\lambda \theta},$$
where $I$ is the identity matrix. The characteristic equation 
\begin{align}\label{eq:TB:ldp:chareq}
det(\Delta(\lambda))=0, \qquad \lambda \in \mathbb{C}
\end{align}
has one solution as zero and all other solutions have negative real parts.
\end{assumpt}

Roughly speaking, after the initial transients have decayed, significant changes in $x$ occur on times of order $O(1/\eps)$ due to the effect of $G$. Since $\sigma$ is mean-zero, large deviations from the corresponding deterministic system are rare on times of order $O(1/\eps)$. We obtain the rate function governing the large deviations.


\section{Spectral theory for DDE}\label{sec:TB:ldp:spectheo}

Under assumption \ref{ass:TB:ldp:assumptondetsys} the space $\C$ can be split as $\C=P\oplus Q$ such that for the unperturbed system $\dot{x}=L_0(\pj_tx)$, the projection of $\pj_tx$ onto $P$ does not change at all, and the norm of the projection of $\pj_tx$ onto $Q$ decays exponentially fast. When the perturbations are present as in \eqref{eq:TB:ldp:main}, the $P$ projection evolves slowly and the $Q$ projection stays small. The $P$ space is one-dimensional. We find a one dimensional evolution equation which is exponential equivalent to the $P$ projection and for which the results from chapter 7 of \cite{FWbook} applies yielding the large deviations rate function.

Here we show, given an $\eta\in \C$, how to find the projection  onto the space $P$. For details, see chapter 7 of \cite{Halebook} and chapter 4 of \cite{Diekmanbook}. 

We use $\R^{n*}$ to distinguish the set of $1\times n$ vectors, from $\R^n$ which is the set of $n\times 1$ vectors. Define the bilinear form  $\la \cdot,\cdot\ra :C([0,r];\R^{n*})\times C([-r,0],\R^n) \to \R$, given by
\begin{equation}\label{eq:TB:ldp:bilinform}
\la \psi,\eta \ra := \psi(0)\eta(0)-\int_{-r}^0\int_0^\theta \psi(s-\theta)d\mu(\theta)\eta(s)ds.
\end{equation}
Choose $\ubar{d}$ such that $\Delta(0)\ubar{d}=0_{n\times 1}$ and $\ubar{d_2}$ such that $\ubar{d_2}\Delta(0)=0_{1\times n}$. Define $\Phi\in \C$ by the constant $\Phi(\bullet)=\ubar{d}$ and $\Psi \in C([0,r];\R^{n*})$ by $\Psi(\bullet)=c\ubar{d_2}$ where the constant $c$ is choosen so that $\la \Psi,\Phi \ra =1$ for the bilinear form in \eqref{eq:TB:ldp:bilinform}. The space $\C$ can be split as $\C=P\oplus Q$ where $P$ is the space spanned by the constant function $\Phi$. The projection operator is $\pi:\C \to P$ given by $\pi(\eta)=\Phi\la \Psi,\eta\ra$. The space $Q$ can be written as $\{\eta \in \C\,:\,\pi(\eta)=0\}$. We find use for $\Psiz \overset{\text{def}}=\Psi(0)$.

The solution to the unperturbed system 
\begin{align}\label{eq:TB:ldp:detDDE}
\dot{x}(t)=L_0(\pj_tx)
\end{align}
can be written as $$\pj_t x=\pi\pj_t x+ (I-\pi)\pj_tx=\Phi z_t + y_t$$ where $z_t=\la \Psi, \pj_t x\ra$ and $y_t=\pj_t x-\Phi z(t)$. Note that $z\in \R$ is a scalar, and $\Phi z_t\in P$ and $y_t\in Q$. It can be shown that for the unperturbed system \eqref{eq:TB:ldp:detDDE}, $\dot{z}=0$, i.e., $z$ is a constant in time. Further, it can be shown that $||y_t||$ decreases to zero exponentially fast (because the dynamics on $Q$ is governed by eigenvalues with negative real parts). Let $\{T(t)\}_{t\geq 0}$, $T(t):\C \to \C$ be the semigroup generated by the DDE \eqref{eq:TB:ldp:detDDE}, i.e. for $\eta\in \C$, $T(t)\eta$ is the solution to \eqref{eq:TB:ldp:detDDE} with the initial condition $\eta$. Then, for $\eta\in P$, $T(t)\eta=\eta$ and, $\exists K,\kappa>0$ such that
\begin{align}\label{eq:TB:ldp:expdecaystatement}
||T(t)\eta||\leq Ke^{-\kappa t}||\eta||, \qquad \qquad \forall \eta\in Q.
\end{align}

Solution to the perturbed equation \eqref{eq:TB:ldp:main} can be written in terms of the semigroup $T$. For this purpose, let $\Ind:[-r,0]\to\R^{n\times n}$ be defined as
$\Ind(\theta)=0_{n\times n}$ for $\theta<0$ and $\Ind(0)=I_{n\times n}$. The solution to \eqref{eq:TB:ldp:main} with initial condition $\pj_0x=\eta$ can be written as
$$\pj_tx(\theta)=T(t)\eta(\theta)+\int_0^tT(t-s)\Ind(\theta) \bigg(G(\pj_sx)+F(\pj_sx)\sigma(\gnoise_s)\bigg)ds, \qquad \theta\in [-r,0].$$
The $j^{th}$ column of $T(t-s)\Ind$ is the solution of \eqref{eq:TB:ldp:detDDE} with the initial condition as the $j^{th}$ column of $\Ind$. Though $(\Ind)_{.j}$ does not belong to $\C$, the bilinear form \eqref{eq:TB:ldp:bilinform} still makes sense and we have $\pi((\Ind)_{.j})=\Phi\Psiz (\Ind)_{.j}$. We still have the exponential decay $||T(t)(I-\pi)(\Ind)_{.j}||\leq Ke^{-\kappa t}||(I-\pi)(\Ind)_{.j}||=Ce^{-\kappa t}$.

Using the fact that $T$ commutes with $\pi$ we have the equations
\begin{align}\label{eq:TB:ldp:Pprojeq}
dz_t\,\,=\,\,{\eps}\Psiz G(\Phi z_t+y_t)dt\,\,+\,\,{\eps}\Psiz F(\Phi z_t + y_t)\sigma(\gnoise_t)dt,
\end{align}
\begin{align*}
y_t = T(t)y_0 + \eps\int_0^tT(t-s)(I-\pi)\Ind\bigg(G(\Phi z_s + y_s)+F(\Phi z_s + y_s)\sigma(\gnoise_s)\bigg)ds.
\end{align*}
Using the exponential decay of $||T(t-s)(I-\pi)\Ind||$ and the boundedness of $F, G, \sigma$ we have that 
\begin{align}\label{eq:TB:ldp:ysmall}
||y_t-T(t)y_0||\,<\, C\eps
\end{align}
for some  $C>0$.


\section{An exponentially equivalent process}

Let the scalar process $\zedx$ be defined by
\begin{align}\label{eq:TB:ldp:Pprojeq_expeq}
d\zedx_t\,\,=\,\,{\eps}\Psiz G(\Phi \zedx_t)dt\,\,+\,\,{\eps}\Psiz F(\Phi \zedx_t)\sigma(\gnoise_t)dt, \qquad \zedx_0=z_0.
\end{align}
Using the bounded derivatives of $F,G$ and boundedness of $\sigma$, and then using \eqref{eq:TB:ldp:ysmall} and the exponential decay \eqref{eq:TB:ldp:expdecaystatement} we have
\begin{align*}
|\zedx_t-z_t|\,\,\leq\,\,C\eps\int_0^t(|\zedx_s-z_s|+||y_s||)ds\,\, &\,\,\leq\,\,C\eps\int_0^t(|\zedx_s-z_s|+||T(s)y_0||+C\eps)ds\\
&\,\,\leq\,\,C\eps^2t + C\eps(1-e^{-\kappa t})+C\eps\int_0^t|\zedx_s-z_s|ds.
\end{align*}
Using Gronwall inequality we have that $\exists C>0$ such that
\begin{align}\label{eq:TB:ldp:Pprojeq_expeq_state}
|z_t-\zedx_t|\leq C\eps, \qquad t\in [0,T/\eps]
\end{align}
for some fixed $T,\eps_0>0$ and all $\eps>\eps_0$. It is easy to see from \eqref{eq:TB:ldp:Pprojeq_expeq} that significant changes for $\zedx$ happens on time of order $O(1/\eps)$, and because $\sigma$ is mean-zero function, significant deviations from the deterministic system $d\zedx_t\,\,=\,\,{\eps}\Psiz G(\Phi \zedx_t)dt$ would be rare on times of order $O(1/\eps)$. By \eqref{eq:TB:ldp:Pprojeq_expeq_state} analogous statement holds for $z_t$.
So we define $z^\eps_t=z_{t/\eps}$ and study the rate function governing the large deviations of $z^\eps$ from the correpsonding determinstic system for $t\in [0,T]$. Define $\zed_t=\zedx_{t/\eps}$. Then, by \eqref{eq:TB:ldp:Pprojeq_expeq_state}, $\zed$ and $z^\eps$ are exponentially equivalent, and so the rate function for $\zed$ and $z^\eps$ are same.

Note that $\zed_t$ is governed by
\begin{align}\label{eq:TB:ldp:Pprojeq_expeq_rescale}
d\zed_t\,\,=\,\,\Psiz G(\Phi \zed_t)dt\,\,+\,\,\Psiz F(\Phi \zed_t)\sigma(\gnoise^\eps_t)dt, \qquad \zed_0=z_0,
\end{align}
where $\gnoise^\eps_t=\gnoise_{t/\eps}$.
The results of Freidlin-Wentzell (chapter 7 of \cite{FWbook}) apply for the large deviations of $\zed_t$ from the deterministic system $\dot{\zedx}_t\,=\,\Psiz G(\Phi \zedx_t)$.


\section{Large deviations of $\zed_t$}

Theorem 7.4.1 in \cite{FWbook} gives the following result.

\begin{thm}
Let the process $\zed$ be governed by \eqref{eq:TB:ldp:Pprojeq_expeq}. Assume the noise $\gnoise$ is homogenous markov process such that for any $\zedx,\alpha\in \R$
\begin{align*}
\lim_{T\to\infty}\frac{1}{T}\ln\expt_{\xi_0}\exp\left(\alpha\int_0^T\Psiz F(\Phi \zedx)\sigma(\gnoise_s)ds\right)=H_F(\zedx,\alpha)
\end{align*}
uniformly in the initial condition $\gnoise_0$ and the function $H_F$ be differentiable with respect to $\alpha$. Let $H(\zedx,\alpha)=\alpha \Psiz G(\Phi \zedx)+H_F(\zedx,\alpha)$. Let $L(\zedx,\beta):=\sup_{\alpha}[\alpha\beta-H(\zedx,\alpha)]$. On $C([0,T];\R)$ introduce the functional
\begin{align*}
S_{0T}(\varphi)=\begin{cases}\int_0^TL(\varphi_s,\dot{\varphi}_s)ds, \qquad \varphi \text{ is absolutely continuous }\\ \infty \qquad \qquad \qquad \text{ otherwise}. \end{cases}
\end{align*}
The functional $S_{0T}$ is the normalized action functional in $C([0,T];\R)$ for the family of processes $\zed$ as $\eps\to 0$, the normalizing coefficient being $1/\eps$. 
\end{thm}

\begin{rmk}
Writing $x^\eps(t)=x(t/\eps)$ we have $x^\eps(t)=\Phi(0) z^\eps_t+y^\eps_t(0)$ and so 
$$|x^\eps(t)-\Phi(0)\zed_t-T(t/\eps)y_0(0)|\,\,\leq\,\,C|z^\eps_t-\zed_t|+C|y^\eps_t-T(t/\eps)y_0|\,\,\leq\,\,C\eps.$$
Recalling that $||T(t/\eps)y_0||\leq Ke^{-\kappa t/\eps}||y_0||$; if $||y_0||$ is small enough, we can approximate the exit rates of $x^\eps$ by exit rates of $\Phi \zed$.
\end{rmk}

\begin{rmk}
For the case of noise $\gnoise$ being $N$-state continuous time Markov chain, theorem 7.4.2 of \cite{FWbook} shows that $H_F(\zedx,\alpha)$ is the largest eigenvalue of the $N\times N$ matrix $Q^{\alpha,\zedx}$ defined by $(Q^{\alpha,\zedx}-Q)_{ij}=\delta_{ij}\sigma_i\alpha \Psiz F(\Phi \zedx)$ where $Q$ is the generator of the Markov chain and $\sigma_i$ is the value of $\sigma$ for the $i^{th}$ state.
\end{rmk}

\begin{rmk}\label{rmk:TB:ldp:instab:2state}
Let $\gnoise$ be a two-state symmetric markov chain with switching rate $g/2$, i.e. 
\begin{align}\label{eq:rateofswitch2state}
\lim_{t\downarrow 0}\frac1t P_{1\to 2}(t)=g/2=\lim_{t\downarrow 0}\frac1t P_{2\to 1}(t)
\end{align}
 where $P_{i\to j}(t)$ is the probability of transition from state $i$ to state $j$ in time $t$. Let $\sigma(\gnoise=1)=-\sigma(\gnoise=2)=\sigma_0$. In this case, the functional $S_{0T}$ can be explicitly evaluated as
\begin{align*}
S_{0T}(\varphi)=\int_0^T\frac{g}{2}\left(1-\sqrt{1-\left(\frac{\dot{\varphi}_s-\Psiz G(\Phi \varphi_s)}{\sigma_0 \Psiz F(\Phi \varphi_s)}\right)^2}\right)ds
\end{align*}
for $\varphi$ absolutely continuous with $|\dot{\varphi}_s-\Psiz G(\Phi \varphi_s)|\leq |\sigma_0 \Psiz F(\Phi \varphi_s)|$ for $s\in [0,T]$ and $\infty$ for all other $\varphi$. The following function would be useful in studying exit related problems:
\begin{align*}
V(t,a,b)=\inf_{\varphi_0=a, \varphi_t=b}S_{0t}(\varphi).
\end{align*}
The solution can be written as
\begin{align*}
V(t,a,b)=\inf_{\substack{\dot{\varphi}_s=\Psiz G(\Phi \varphi_s)+\sigma_0 \Psiz F(\Phi \varphi_s)u_s,\\ |u_s|\leq 1,  \,\,\varphi_0=a,\,\, \varphi_t=b}}\int_0^t\frac{g}{2}\left(1-\sqrt{1-u^2_s}\right)ds.
\end{align*}
\end{rmk}


\section{Linear delay equations with fast markov perturbations}

In this section we make an independent observation regarding processes of the form
\begin{align}\label{eq:TB:ldp:main:fastSwitch}
\dot{x}^\eps(t)=L_0(\pj_tx^\eps)+\sigma(\gnoise^\eps_t), \qquad \quad \pj_0x^\eps=\eta\in \C
\end{align}
where $\gnoise^\eps_t=\gnoise_{t/\eps}$ with $\gnoise$ being a homogenous markov process and $\sigma$ being a mean-zero $\R^n$-valued function of the noise. Assume that for any $\alpha\in \R^n$
\begin{align*}
\lim_{T\to\infty}\frac{1}{T}\ln\expt_{\xi_0}\exp\left(\int_0^T\alpha^*\sigma(\gnoise_s)ds\right)=H(\alpha)
\end{align*}
uniformly in the initial condition $\gnoise_0$ and the function $H$ be differentiable with respect to $\alpha$.  Let $L(\beta):=\sup_{\alpha}[\alpha^*\beta-H(\alpha)]$. On $C([0,T];\R^n)$ introduce the functional
\begin{align*}
S^{\sigma}_{0T}(\varphi)=\begin{cases}\int_0^TL(\dot{\varphi}_s)ds, \qquad \varphi \text{ is absolutely continuous }\\ \infty \qquad \qquad \qquad \text{ otherwise}. \end{cases}
\end{align*}
The functional $S^{\sigma}_{0T}$ is the normalized action functional in $C([0,T];\R^n)$ for the family of processes $\int_0^{\cdot}\sigma(\xi^\eps_s)ds$ as $\eps\to 0$, the normalizing coefficient being $1/\eps$. 

Define the map $\mathfrak{B}_{\eta}:C([0,T];\R^n)\to C([0,T];\R^n)$ by $\mathfrak{B}_{\eta}\psi=v$ where $v$ is the solution of 
\begin{align*}
v(t)=\eta(0)+\int_0^tL_0(\pj_sv)ds+\psi(t).
\end{align*}
with the understanding that $\pj_0v=\eta$.
More explicit representation of $v$ can be given by the variation-of-constants formula. The map $\mathfrak{B}_{\eta}$ has inverse given by $(\mathfrak{B}_{\eta}^{-1}v)(t)=v(t)-\eta(0)-\int_0^tL_0(\pj_sv)ds$. It can be shown using Gronwall inequality that $\mathfrak{B}_{\eta}$ is Lipschitz. By contraction principle we have that the action functional for $x^\eps$ is given by
\begin{align*}
S_{0T}(\varphi)=\begin{cases}\int_0^TL\left(\dot{\varphi}_s-L_0(\pj_s\varphi)\right)ds, \qquad \varphi \text{ is absolutely continuous }\\ \infty \qquad \qquad \qquad \text{ otherwise}, \end{cases}
\end{align*}
with the understanding that $\pj_0\varphi$ is the initial condition.

Consider the case of $x$ being $\R$-valued, and $\xi$ being a two-state markov chain as in the remark \ref{rmk:TB:ldp:instab:2state}. The following function would be useful in studying exit related problems:
\begin{align*}
V(t,\eta,b)=\inf_{\pj_0\varphi=\eta, \,\,\varphi(t)=b}S_{0t}(\varphi).
\end{align*}
The solution can be written as
\begin{align*}
V(t,\eta,b)=\inf_{\substack{\dot{\varphi}_s=L_0(\pj_s\varphi)+\sigma_0 u_s,\\ |u_s|\leq 1,  \,\,\pj_0\varphi=\eta,\,\, \varphi(t)=b}}\int_0^t\frac{g}{2}\left(1-\sqrt{1-u^2_s}\right)ds.
\end{align*}
Let $f:[-r,\infty)\to\R$ be defined by $f(t)=0$ for $t<0$, $f(0)=1$, and for $t>0$, $f$ satisfies $\dot{f}(t)=L_0(\pj_tf)$. Let $\{T(t)\}_{t\geq 0}$ be the solution semigroup as defined in section \ref{sec:TB:ldp:spectheo}. 
Then the solution to $$\dot{\varphi}_s=L_0(\pj_s\varphi)+\sigma_0 u_s, \qquad \pj_0\varphi=\eta,$$
can be represented using the variation-of-constants formula as
$$\varphi(t)=T(t)\eta(0)+\int_0^tf(t-s)\sigma_0u_sds.$$
Hence we have
\begin{align*}
V(t,\eta,b)=\inf_{\substack{\int_0^tf(t-s)\sigma_0u_sds\,\,=\,\,b-T(t)\eta(0) \\ |u_s|\leq 1}}\int_0^t\frac{g}{2}\left(1-\sqrt{1-u^2_s}\right)ds.
\end{align*}
The RHS above can be computed explicity using calculus of variations. We have for the optimality, $u_s=\frac{-\rho f(t-s)}{\sqrt{1+\rho^2f^2(t-s)}}$ with the Lagrange multiplier $\rho$ obtained using $\int_0^tf(t-s)\sigma_0u_sds=b-T(t)\eta(0)$.

Note that $\frac{1}{\eps}\int_0^\cdot \sigma(\xi^\eps_s)ds$ converges weakly as $\eps \to 0$ to $\sigma_0\sqrt{g}\,W_{\cdot}$ where $W$ is a Wiener process. However, the large deviations principle for $\dot{x}^\eps(t)=L_0(\pj_tx^\eps)+\eps(\frac{1}{\eps}\sigma(\gnoise^\eps_t))$ is different from the large deviations principle for $dx^\eps(t)=L_0(\pj_tx^\eps)dt+\eps \sigma_0\sqrt{g}dW_t$.

Large deviations for DDE with noise as Wiener process is considered in \cite{SEAMzhangLDPDelay}.

\end{document}